\documentclass[12pt]{amsart}
\usepackage{amssymb}

\usepackage[all]{xy}
\newtheorem{theorem}{Theorem}[section]
\newtheorem{lemma}[theorem]{Lemma}
\newtheorem{proposition}[theorem]{Proposition}
\newtheorem{corollary}[theorem]{Corollary}

\theoremstyle{definition}
\newtheorem{definition}[theorem]{Definition}

\theoremstyle{remark}
\newtheorem{remark}[theorem]{Remark}

\numberwithin{equation}{section}

\begin{document}

\title[Disjointly weak compactness in Banach lattices ]
 {Disjointly weak compactness in Banach lattices}

 \author[B. Xiang]
{Bo Xiang}
\address{School of Mathematics, Southwest Jiaotong
University, Chengdu 610031,  China}
\email{2734598188@qq.com}

\author[J.X.  Chen]
{Jin Xi Chen }

\address{School of Mathematics, Southwest Jiaotong
University, Chengdu 610031, China}
 \email{jinxichen@swjtu.edu.cn}

\author[L. Li]
{Lei Li}
\address{School of Mathematical Sciences and LPMC, Nankai University, Tianjin 300071,  China}
\email{leilee@nankai.edu.cn}
\thanks{ The third author was partly supported  by NSFC (No.12171251).}

%    General info
\subjclass[2010]{Primary 46B42; Secondary 46B50, 47B65}
%\dedicatory{}

\keywords{disjointly weakly compact set, disjointly weakly compact operator, $(L-,\, M-)$\,weakly compact operator, positive operator, Banach lattice}

\begin{abstract}
We give some characterizations of disjointly weakly compact sets in Banach lattices, namely, those sets in whose solid hulls every disjoint sequence   converges weakly to zero. As an application, we prove that a bounded linear operator from a Banach space to a Banach lattice  is an almost $(L)$ limited operator if and only if it is a disjointly weakly compact operator, indeed, an operator which carries bounded sets to disjointly weakly compact ones. Some results on weak precompactness and ($L$\,-, $M$-)weak compactness of  disjointly weakly compact operators are also obtained.
	
\end{abstract}

\maketitle \baselineskip 4.95mm

\section{Introduction}

\par Throughout this paper, $X, Y$ will denote real Banach spaces and $E,F$  will denote real Banach lattices. We denote by $B_{X}$ the closed unit ball of $X$. $Sol(A)$ denotes the solid hull of a subset $A$ of a Banach lattice. The positive cone of $E$ will be denoted by $E^{\,+}$.

\par Let us recall that a bounded subset $A$ of $X$ is called a \textit{Dunford-Pettis}  (resp. \textit{limited}) \textit{set} in $X$ if every weakly (resp. weak$^{*}$) null sequence in $X^{\prime}$  converges to zero uniformly on $A$ (see \cite{Andrews, BD}). $X$ has the \textit{Dunford-Pettis property } (resp. \textit{Dunford-Pettis$^{*}$ property}) whenever every relatively weakly compact set in $X$ is a Dunford-Pettis (resp. limited) set. The Dunford-Pettis$^{*}$ property, introduced first by Borwein, Fabian and Vanderwerff \cite{BFV}, is stronger than the Dunford-Pettis property. Carri\'{o}n, Galindo and Louren\c{c}o \cite{CGL} showed that $X$ has the Dunford-Pettis$^{*}$ property if and only if every continuous linear operator $T:X\rightarrow c_{0}$ is a Dunford-Pettis operator.

\par The disjoint versions of Dunford-Pettis sets and limited sets in  Banach lattices were introduced respectively by Bouras \cite{Bour} and Chen et al. \cite{CCJ}.  Recall that a bounded subset $A$ of a Banach lattice $E$ is called an \textit{almost Dunford-Pettis set} (resp. \textit{almost limited set}) if every disjoint weakly (resp. weak$^{*}$) null sequence $(f_{n})$ of $ E^{\prime}$ converges uniformly to zero on $A$, that is, $\sup_{x\in A}|f_{n}(x)|\rightarrow0$. A Banach lattice $E$ has the \textit{ weak Dunford-Pettis property } (resp. \textit{weak Dunford-Pettis$^{*}$ property}) if every relatively weakly compact set in $E$ is  almost Dunford-Pettis (resp. almost limited). The dual versions of limited sets and almost limited sets are the so-called $L$\,-sets and almost $L$\,-sets, respectively.  A bounded subset $B$ of the dual $X^{\prime}$ of a Banach space $X$ is called an  \textit{$L$\,-set} in $X^{\prime}$ whenever every  weakly null sequence $(x_{n})$ of $X$ converges uniformly to zero on the set $B$  (see \cite{Em}).  A norm bounded subset $B$ of the dual $E^{\prime}$ of a Banach lattice $E$ is called an  \textit{almost $L$\,-set} in $E^{\prime}$ whenever every  disjoint weakly null sequence $(x_{n})$ of $E$ converges uniformly to zero on $B$ (see \cite{Aqzzouz1}). For simplicity, a sequence $(f_{n})$ of continuous functionals is called an  $L$\,-sequence if the set $\left\{{{f_n}:{n}\in{\mathbb{N}}}\right\}$ is an $L$\,-set.
\par Subsequently, some special classes of operators related to Dunford-Pettis sets, limited sets and  $L$\,-sets appeared successively in the literature. A bounded linear operator $T:X\to{Y}$ is called
 \begin{itemize}
   \item[$\raisebox{0mm}{------}$]  a \textit{limited operator} if $f_{n}\xrightarrow {w^{*}}0$ in $Y^{\prime}$ implies $ \lim\|T^{\prime}f_{n}\|=0$ or, equivalently, if $T(B_{X})$ is a limited set in $Y$ (see \cite{BD}, the dual version of Dunford-Pettis operators).
   \item[$\raisebox{0mm}{------}$] a \textit{Dunford-Pettis completely continuous} (\,DPcc for short) \textit{operator} if it maps Dunford-Pettis and weakly null sequences to norm null ones (see \cite{WC}).

\end{itemize}
Also, a bounded linear operator $T:E\to{X}$ from a Banach lattice into a Banach space is called
\begin{itemize}

	\item[$\raisebox{0mm}{------}$]  \textit{almost Dunford-Pettis}, if it maps disjoint and weakly null sequences to norm null ones (see \cite{Sa}, the disjoint version of Dunford-Pettis operators)

	\item[$\raisebox{0mm}{------}$]   \textit{almost Dunford-Pettis completely continuous} (\,aDPcc for short),  if it maps  Dunford-Pettis and weakly null sequences consisting of disjoint terms to norm null ones ( see \cite{FMM}, the disjoint version of DPcc operators).
\end{itemize}
 An operator $T:X\to{E}$ is called \textit{almost limited} if $T(B_{X})$ is an almost limited set in $E$, or equivalently, if $\|{T^{\prime}f_{n}}\|\to{0}$ for every disjoint weak$^*$ null sequence $(f_n)$ in $E^{\prime}$ (see \cite{El, MEM}). Clearly, this class of operators is  the disjoint version of limited operators.

\par Recently,  Oughajji, El Fahri and Moussa \cite{Oughajji} considered  the dual version of aDPcc operators which is at the same time a weak version of almost limited operators to some extent, and introduced the class of almost (\textit{L}) limited operators. An operator $T$ from a Banach space $X$ into a Banach lattice $E$ is said to be \textit{almost (L) limited} if $\|T^{\prime}f_n\|\rightarrow0$  for every  disjoint
$L$\,-sequence $(f_n)$ in $E^{\prime}$ satisfying $f_{n}\xrightarrow {w^{*}}0$. By definition, in order to know about an almost $(L)$ limited operator the authors in \cite{Oughajji} had to study the behavior of its adjoint. There, some interesting results were obtained when the spaces in question are $\sigma$-Dedekind complete. Among other things, they proved that every disjoint sequence $(x_n)\subset Sol(T(B_E))$
is weakly null if $T:E\to{F}$ is an almost $(L)$ limited operator between two Banach lattices. Following W. Wnuk \cite{Wnuk} we call a bounded subset $A$  of a Banach lattice  \textit{disjointly weakly compact} if every disjoint sequence $(x_n)\subset Sol(A)$ is weakly null. Now, in terms of disjointly weakly compact sets  the authors of \cite{Oughajji} indeed showed that if $T:E\to{F}$ is an almost $(L)$ limited operator then $T(B_E)$ is disjointly weakly compact. We conjecture that the converse also holds. To establish the conjecture,  we first give some characterizations of disjointly weakly compact sets in Banach lattices. Then, naturally, we introduce the class of disjointly weakly compact operators, namely, operators which carry bounded sets to disjointly weakly compact ones. We prove that a bounded linear operator $T:X\to{E}$  is almost $(L)$ limited if and only if it is disjointly weakly compact.  Some results on weak precompactness and ($L$\,-, $M$-) weak compactness of  disjointly weakly compact operators are also obtained.

\par Our notions are standard. For the theory of Banach lattices and  operators, we refer the reader to the monographs \cite{AB, Meyer, Wnuk}.

\section{Disjointly weak compactness in Banach lattices}

\par Let us recall that a norm bounded subset $A$ of a Banach lattice $E$ is said to be \textit{$L$\,-weakly compact} if $\|{x_{n}}\|\to{0}$ for every disjoint sequence $(x_{n})$ in $Sol(A)$. Here the letter ``$L$" comes from the term $L$\,-space because $L$\,-weakly compact sets are automatically relatively weakly compact, and  relatively weakly compact sets and $L$\,-weakly compacts sets coincide in $L$\,-spaces (see, e.g., \cite{Meyer}).  In \cite{Wnuk}  Wnuk introduced a special class of sets  that he called ``disjointly weakly compact sets", which is obviously the weak version of $L$\,-weakly compact sets.
\begin{definition}\cite{Wnuk}\label{Definition 2.1}
	A norm bounded subset $A$ of a Banach lattice $E$ is called disjointly weakly compact if $x_{n} \stackrel{w}{\rightarrow}0$ for disjoint sequences $(x_{n})\subset{Sol(A)}$.	
\end{definition}
Although afterwards the name  of \textit{disjointly weakly compact set} seldom appeared in the literature, we know that the property described in Definition \ref{Definition 2.1} is shared by many other classes of sets. For example, every order bounded set and every weakly sequentially precompact set (and hence relatively weakly compact set) in a Banach lattice is disjointly weakly compact (see, e.g., \cite[Proposition 2.5.12 iii)]{Meyer}). It is well known that the closed unit ball $B_E$ of a Banach lattice $E$ is disjointly weakly compact if and only if $E^{\prime}$ has order continuous norm (see, e.g., \cite[Theorem 2.4.14]{Meyer}).
\par For our further discussions we need the following lemma due to the authors \cite[Lemma 2.2]{Xiang}.
\begin{lemma}\label{Lemma 2.2}\cite{Xiang} Let $E$ be a Banach lattice and let $A$ be an almost $L$\,-set in $E^{\prime}$.
	\begin{enumerate}
		\item If $(f_n)$ is a sequence in $Sol(A)$ satisfying ${\vert{f_n}\vert}\stackrel{w^*}{\longrightarrow}0$ , then $(f_n)$ is  an $L$-sequence.
		\item If $(f_n)$ is a disjoint sequence in Sol$(A)$, then $\vert{f_n}\vert\stackrel{w^*}{\longrightarrow}0$ and $(f_n)$ is  an $L$-sequence.
	\end{enumerate}
\end{lemma}
  By this lemma every disjoint $L$\,-sequence in $E^{\prime}$ is automatically (absolutely) weak$^*$ null. Also, every order bounded disjoint sequence is an $L$\,-sequence. Now we turn our attention to the characterizations of disjointly weakly compact sets.
\begin{theorem}\label{Theorem 2.3}
	Let $A$ be a norm bounded subset of a Banach lattice $E$. The following statements are equivalent.
	\begin{enumerate}
		\item $A$ is disjointly weakly compact.
		\item Every order bounded sequence $(f_{n})\subset{(E^{\prime})^{+}}$ satisfying $f_{n} \stackrel{w^*}{\rightarrow}0$ converges uniformly to zero on $A$.
		\item Every disjoint $L$\,-sequence $(f_{n})\subset{E^{\prime}}$ converges uniformly to zero on $A$.
		\item Every order bounded disjoint sequence $(f_{n})\subset{E^{\prime}}$ converges uniformly to zero on $A$.

	\end{enumerate}	
\end{theorem}
\begin{proof}
	(1)$\Rightarrow$(2)
	Let $(f_{n})_{n=1}^{\infty}$ be an order bounded weak$^*$ null sequence in ${(E^{\prime})^{+}}$ and let $0\leq f_{n}\leq f$ for some $f\in (E^{\prime})^{+}$.  Let $A$ be a disjointly weakly compact subset of $E$.  By definition, $Sol(A)$ is a disjointly weakly compact subset in $E$. Now we can assume without loss of generality that $A$ is solid. To finish the proof, we have to show that $\rho_{A}(f_{n}):=\sup_{x\in{A}}{|f_{n}(x)|}\to{0}$.  It follows from Proposition 2.3.4  of \cite{Meyer} that there exists a disjoint sequence $(x_{n})_{n=1}^{\infty}\subset{A\cap{E^{+}}}$ such that
	$$\limsup_{n\rightarrow\infty}\rho_{A}(f_{n})=\limsup_{n\rightarrow\infty}\langle f_{n}, x_{n}\rangle\le\lim_{n\rightarrow\infty}{f(x_{n})}=0$$  since $A$ is disjointly weakly compact and ${x_{n}} \stackrel{w}{\rightarrow}0$.

\par (2)$\Rightarrow$(4) Let $(f_{n})$ be an order bounded disjoint sequence of $E^{\prime}$. Since both $(f_{n}^{+})$ and $(f_{n}^{-})$ are disjoint and order bounded  in ${(E^{\prime})^{+}}$. Hence, $(f_{n}^{+})$ and $(f_{n}^{-})$ are both weakly null. Then from the inequality
	$$\begin{aligned}
		\sup_{x\in{A}}{\vert{f_{n}(x)}\vert}\le\sup_{x\in{A}}{\vert{f_{n}^{+}(x)}\vert}+\sup_{x\in{A}}{\vert{f_{n}^{-}(x)}\vert}
	\end{aligned}$$
	it follows that $\sup_{x\in{A}}{\vert{f_{n}(x)}\vert}\to{0}$.
\par (1)$\Leftrightarrow$(4) For every $0\leq f\in E^{\prime} $ let $B=[-f, f\,]$. Then the equivalence follows easily from a well known result due to O. Burkinshaw and P. G. Dodds (see, e.g., \cite[Theorem 5.63]{AB}\,).

\par(1)$\Rightarrow$(3) Let $(f_{n})_{n=1}^{\infty}$ be an arbitrary disjoint $L$\,-sequence in $E^{\prime}$. We have to show that $(f_{n})$ converges uniformly to zero on $A$. Otherwise, passing to a subsequence if necessary, there would  some $\epsilon>{0}$ such that $\sup_{x\in{A}}{\vert{f_{n}(x)}\vert}>{\epsilon}$ for all $n$. Therefore, there exists a sequence $(x_{n})_{n=1}^{\infty}$ in $A$ such that $\vert{f_{n}}\vert(\vert{x_{n}}\vert)\geq\vert{f_{n}}({x_{n}})\vert>{\epsilon}$ for every $n$. Since ${\{f_n:{n}\in{\mathbb{N}}\}}$ is a disjoint $L$\,-sequence, it follows from Lemma \ref{Lemma 2.2}\,(2) that $\vert{f_{n}}\vert\stackrel{w^*}{\longrightarrow}0$.  Now, as we did in the proof of  Lemma 2.2 of \cite{Xiang}  we can find a strictly increasing subsequence $(n_{k})_{k=1}^{\infty}$ of $\mathbb{N}$  and a disjoint sequence $(y_{m})$ in ${Sol(A)}\cap E^+$ such that
$$0\leq y_m \leq |x_{n_{m+1}}| \quad \textmd{and} \quad |f_{n_{m+1}}|(y_m)>\frac{\varepsilon}{2} \,\,( m \,\,\textmd{sufficiently\, large})$$By the well-known Riesz-Kantorovich Formula pick a sequence $(u_m)$ in $E$ such that  ${\vert{u_m}\vert}\le{y_m}$ and ${\vert{f_{n_{m+1}}(u_m)}\vert}>{\frac{\epsilon}{2}}$ for $m$ sufficiently large. On the other hand, it is clear that $(u_m)\subset{Sol(A)}$ and $(u_m)$ is also disjoint. It follows that $u_{m}\stackrel{w}{\rightarrow}0$ since $A$ is disjointly weakly compact.  Because $(f_{n})_{n=1}^{\infty}$ is an $L$\,-sequence in $E^{\prime}$, by definition $$\vert{f_{n_{m+1}}(u_m)}\vert\leq \sup_{n}\vert f_{n}(u_m)\vert\to{0}\qquad (m\rightarrow\infty)$$This leads to a contradiction.

\par (3)$\Rightarrow$(4) This follows easily because every order bounded disjoint sequence of ${E^{\prime}}$ is an $L$\,-sequence.
\end{proof}

\begin{remark}\label{Remark 2.4}
\par(1) Clearly, every almost Dunford-Pettis set (and hence every almost limited set) is  disjointly weakly compact. This can be easily seen  from Theorem \ref{Theorem 2.3} (1)$\Leftrightarrow$(4). The converse does not hold in general. For instance,  the closed unit ball $B_{\ell_{2}}$ is weakly compact, hence disjointly weakly compact, but it is not almost Dunford-Pettis.
	\par (2) Let us recall that a bounded subset $A$ of a Banach space $X$ is said to be weakly (sequentially) precompact provided  that every sequence from $A$ has a weakly Cauchy subsequence.     Ghenciu \cite{Ghenciu} gave a  characterization of weakly  precompact sets in Banach spaces. He proved that  a bounded subset $A$ of $X$ is weakly precompact  if and only if each weak$^{*}$ null $L$\,-sequence $(f_n)$ in $X^{\prime}$ converges to zero uniformly on $A$. From this we can see that every weakly precompact set in a Banach lattice is  disjointly weakly compact. This is also a known result included in   Proposition 2.5.12 iii) of \cite{Meyer}. The converse does not hold in general. For example, $B_{\ell_{\infty}}$ is  disjointly weakly compact in $\ell_{\infty}$. However, by Rosenthal's $\ell_{1}$-theorem $B_{\ell_{\infty}}$ is not weakly precompact. Generally, every disjointly weakly compact set in a Banach lattice $E$  is weakly precompact if and only if  every order interval in $E$ is weakly precompact. This result is due to the authors \cite{Xiang}.
\end{remark}

In terms of weakly Cauchy sequences  we characterize disjointly weakly compact sets  as follows.
\begin{theorem} Let $A$ be a bounded subset of a Banach lattice $E$. $A$  is disjointly weakly compact if and only if  $(x_{n})$ is weakly Cauchy for every disjoint sequence $(x_{n})\subset{Sol(A)}$.
\end{theorem}
\begin{proof}
  We only have to verify the sufficiency. If every disjoint sequence from ${Sol(A)}$ is weakly Cauchy, we assume by way of contradiction that  $A$ fails to be disjointly weakly compact. Then, as we did in the proof of Corollary 2.5 of \cite{Xiang}, we could find a disjoint sequence $(x_{n})_{1}^{\infty}\subset Sol(A)^{+}$ which is equivalent to the unit vector basis $(e_{n})_{1}^{\infty}$ of $\ell_{1}$. Therefore, the disjoint sequence $(x_{n})_{n}^{\infty}$ is not weakly Cauchy. This leads to a contradiction.
\end{proof}

Note that every relatively weakly compact set in a Banach lattice is disjointly weakly compact. The following result gives the condition under which the class of relatively weakly compact sets and the class of disjointly weakly compact sets coincide.
\begin{proposition}\label{Proposition 2.5}
	Let $E$ be a Banach lattice. Every disjointly weakly compact set in $E$ is relatively weakly compact if and only if $E$ is a $KB$-space.
\end{proposition}
\begin{proof} If every disjointly weakly compact set in $E$ is relatively weakly compact, then every  weakly precompact set in $E$ is relatively weakly compact. Therefore, $E$ is weakly sequentially complete, that is, $E$ is a $KB$-space (see \cite[Theorem 4.60 ]{AB}). For the converse, let $E$ be a $KB$-space and let $A$ be a disjointly weakly compact set in $E$. We can assume without loss of generality that $A$ is solid. It follows from Theorem \ref{Theorem 2.3} (4) that $\rho_{A}(f_{n})\to{0}$ for every order bounded disjoint sequence $(f_{n})\subset{E^{\prime}}$. Thus the desired  result follows from Theorem 2.5.3 of \cite{Meyer}.
\end{proof}	
\par
 Recall that a Banach lattice $E$ has the \textit{positive Schur property} or \textit{weak Schur property} if every weakly null sequence with positive terms is norm null, equivalently, if every disjoint weakly null sequence in $E$ is norm null. It should be noted  that every Banach lattice with the positive Schur property is a $KB$-space.  By \cite[Corollary 3.6.8]{Meyer} $E$ has the positive Schur property if and only if every relatively weakly compact subset of $E$ is $L$-weakly compact. The following is an immediate consequence of Proposition \ref{Proposition 2.5}.
 \begin{corollary}\label{psp}
  A Banach lattice $E$ has the positive Schur property if and only if every disjointly weakly compact subset of $E$ is $L$-weakly compact.
\end{corollary}
\par Let us recall from \cite{CCJ} that a Banach lattice $E$ has the weak Dunford-Pettis$^*$ property (wDP$^*$ property for short) if every relatively weakly compact set in $E$ is almost limited, or equivalently, if $f_{n}(x_{n})\to{0}$ for every weakly null sequence $(x_{n})$ in $E$ and for every disjoint weak$^*$ null sequence $(f_{n})$ in $E^{\prime}$.
\begin{theorem}\label{Theorem 2.6}
	A Banach lattice $E$ has the wDP$^{\,*}$ property if and only if every disjointly weakly compact subset of $E$ is almost limited.
\end{theorem}
\begin{proof}It suffices to prove the  ``only if " part. Let $A$ be a disjointly weakly compact subset of $E$ and let $(f_{n})$ be an arbitrary disjoint weak$^*$ null sequence in $E^{\prime}$.  Since $E$  has the wDP$^*$ property, then for every weakly null sequence $(x_{n})$ in $E$ we have ${f_{n}(x_{n})}\to{0}$, and hence $(f_{n})$ is an $L$\,-sequence (cf. \cite[Proposition 2.2]{Aqzzouz1}). Since $A$ is  disjointly weakly compact, it follows from Theorem \ref{Theorem 2.3} (3) that ${\sup_{x\in{A}}{\vert{f_{n}(x)}\vert}}\to{0}$. This implies that $A$ is almost limited.
\end{proof}	

\par We say that a Banach lattice $E$ has the weak Dunford-Pettis property (wDP property for short) if every weakly compact operator from $E$ into any Banach space $Y$ is almost Dunford-Pettis (see \cite{Leung}). Later,  Bouras \cite{Bour} proved that  $E$ has the weak Dunford-Pettis property (wDP property for short) if each relatively weakly compact set in $E$ is almost Dunford-Pettis, or equivalently, if $f_{n}(x_{n})\to{0}$ holds for each weakly null sequence $(x_{n})$ in $E$ and for each disjoint weakly null sequence $(f_{n})$ in $E^{\prime}$.
\begin{theorem}\label{Theorem 2.7}
	For a Banach lattice $E$ the following statements are equivalent.
	\begin{enumerate}
		\item  $E$ has the wDP property.
		\item Each disjointly weakly compact subset of $E$ is an almost Dunford-Pettis set.	
		\item  Each disjointly weakly compact subset of $E^{\prime}$ is  an almost $L$\,-set.
	\end{enumerate}	
\end{theorem}
\begin{proof}
	(1)$\Leftrightarrow$(2) The proof is very similar to  that of Theorem \ref{Theorem 2.6}.
\par
	(3)$\Rightarrow$(1)  Every relatively weakly compact subset of $E^{\prime}$ is almost $L$\,-set since  every disjointly weakly compact set is an almost $L$\,-set. The desired result follows from Theorem 3.10 of \cite{Aqzzouz1}. \par
	(1)$\Rightarrow$(3) Let $A$ be a disjointly weakly compact  subset of $E^{\prime}$. We have to prove that $A$ is an almost $L$\,-set. To this end, let $(x_{n})$ be an arbitrary disjoint weakly null sequence of $E$. If $j:E\to{E^{\prime\prime}}$ is the canonical embedding, then $(jx_{n})$ is also disjoint weakly null sequence of $E^{\prime\prime}$. Since $E$ has the wDP property, then $(jx_{n})(f_{n})=f_{n}(x_{n})\to{0}$ for every weakly null sequence $(f_{n})$ in $E^{\prime}$. Therefore, $(jx_{n})$ is an $L$\,-sequence in $E^{\prime\prime}$. Since $A$ is disjointly weakly compact, we have
	$$\sup\limits_{f\in{A}}{\vert{f(x_{n})}\vert}={\sup\limits_{f\in{A}}{\vert{(jx_{n})(f)}\vert}\to{0}}.$$
This implies that $A$ is an almost $L$\,-set.
\end{proof}

\section{ Disjointly weakly compact operators on Banach lattices}

\par Recall that a bounded linear operator $T:X\to{E}$ from a Banach space into a Banach lattice is called $L$\,-weakly compact if $T(B_{X})$ is an $L$\,-weakly compact subset of $E$, that is, every disjoint sequence in $Sol(T(B_{X}))$ is norm null. Inspired by this we naturally introduce the following definition of disjointly weakly compact operators.
\begin{definition}\label{Definition 3.1}
	A bounded linear operator $T:X\to{E}$ from a Banach space $X$ to a Banach lattice $E$ is said to be disjointly weakly compact if
	$T(B_{X})$ is disjointly weakly compact  in $E$.
\end{definition}
\begin{remark}

  (1) Clearly, every weakly precompact operator and every almost limited operator from a Banach space to a Banach lattice is disjointly weakly compact. If $E^{\prime}$ has order continuous norm, then every bounded linear operator from an arbitrary Banach space $X$ to $E$ is disjointly weakly compact.

   \par (2) Let $S,T:E\to{F}$ be two positive operators between Banach lattices such that $0\le{S}\le{T}$ holds. If $T$ is a disjointly weakly compact operator, then $S$ is likewise disjointly weakly compact. This follows from the fact that $S(B_{E})\subset{Sol(T(B_{E}))}$.

\end{remark}
\par Recently,  Oughajji, El Fahri and Moussa in \cite{Oughajji} introduced the class of  almost (\textit{L}) limited operators.  An operator $T$ from a Banach space $X$ into a Banach lattice $E$ is said to be \textit{almost (L) limited} if $\|T^{\prime}f_n\|\rightarrow0$  for every  disjoint
weak$^*$ null $L$\,-sequence $(f_n)$ in $E^{\prime}$. The authors in \cite{Oughajji} worked under the condition that the spaces involved are assumed to be $\sigma$-Dedekind complete.  They proved that every disjoint sequence $(x_n)\subset Sol(T(B_E))$
is weakly null if $T:E\to{F}$ is an almost $(L)$ limited operator between two Banach lattices. Now, in terms of disjointly weakly compact operators  they indeed showed that if $T:E\to{F}$ is an almost $(L)$ limited operator then $T$ is a disjointly weakly compact operator. Note that every disjoint $L$\,-sequence in $E^{\prime}$ is automatically  weak$^*$ null. The following result shows that a bounded linear operator $T$ from a Banach space to a Banach lattice is an almost $(L)$ limited operator if and only if $T$ is  disjointly weakly compact, or equivalently, if and only if $T^{\prime}$ is order weakly compact.

\begin{theorem}\label{Theorem 3.2}
	Let $T:X\to{E}$ be a bounded linear operator from a Banach space $X$ to a Banach lattice $E$. Then the following statements are equivalent.
	\begin{enumerate}
		\item $T$ is disjointly weakly compact.
       \item $\|{T^{\prime}f_{n}}\|\to{0}$ for each order bounded sequence $(f_{n})\subset{(E^{\prime})^{+}}$ satisfying $f_{n} \stackrel{w^{*}}{\rightarrow}0$.
		\item $\|{T^{\prime}f_{n}}\|\to{0}$ for each disjoint $L$-sequence $(f_{n})\subset{E^{\prime}}$ (hence $f_{n} \stackrel{w^{*}}{\longrightarrow}0$), i.e., $T$ is almost $(L)$ limited.
		\item $\|{T^{\prime}f_{n}}\|\to{0}$ for each order bounded disjoint sequence $(f_{n})\subset{E^{\prime}}$, i.e., $T^{\prime}$ is order weakly compact.
       	\item $\|{T^{\prime}f_{n}}\|\to{0}$ for each weak$^{*}$ null $L$-sequence $(f_{n})\subset{(E^{\prime})^{+}}$.

\end{enumerate}	
\end{theorem}
\begin{proof}
	
	(1)$\Leftrightarrow$(2) $\Leftrightarrow$(3)$\Leftrightarrow$(4) follows from Theorem \ref{Theorem 2.3} and the equalities: $$\|{T^{\prime}f_{n}}\|= {\sup_{x\in{B_{X}}}{|{{T^{\prime}}f_{n}(x)}|}}={\sup_{z\in{T(B_{X})}}{|{f_{n}(z)}|}}.$$\par

	(3)$\Rightarrow$(5) Let $(f_{n})$ be an $L$\,-sequence in $(E^{\prime})^{+}$ satisfying $f_{n} \stackrel{w^{*}}{\longrightarrow}0$, and let $A$ be the solid hull of the set ${\{f_{n}:{n}\in{\mathbb{N}}\}}$. By Lemma \ref{Lemma 2.2}\,(2), every disjoint sequence $(g_{n})\subset{A}$ is an $L$\,-sequence. Then, by hypothesis, we have $\|{T^{\prime}g_{n}}\|\to{0}$. Therefore, it follows from Theorem 4.36 of \cite{AB} that for an arbitrary $\epsilon>0$  there exists some $h\in{(E^{\prime})^{+}}$ lying in the ideal generated by $A$ such that
	$\|{T^{\prime}(f_{n}-h)^{+}}\|<{\epsilon}$ for all $n$. On the other hand, it is clear that $(f_{n}\land{h})$ is an order bounded sequence of ${(E^{\prime})^{+}}$ satisfying ${f_{n}\land{h}}\stackrel{w^*}{\longrightarrow}0$.  Then  the statement (2) implies that $\|{T^{\prime}(f_{n}\land{h})}\|\to{0}$. Now, from the equality $f_{n}=(f_{n}-h)^{+}+(f_{n}\land{h})$ it follows that
	$$\|{T^{\prime}f_{n}}\|\leq{\|{T^{\prime}(f_{n}-h)^{+}}\|+{\|{T^{\prime}(f_{n}\land{h})}\|}}$$which implies that ${\lim\sup_{n\rightarrow\infty}{\|{{T^{\prime}}f_{n}}\|}}\leq{\epsilon}$. Since $\epsilon>{0}$ is arbitrary, we have  $\|{T^{\prime}f_{n}}\|\to{0}$, as desired.\par
	(5)$\Rightarrow$(3) Let $(f_{n})$ be a disjoint $L$\,-sequence in $(E^{\prime})$. By Lemma 2.2\,(2), $(f_{n}^{+})$ and $(f_{n}^{-})$ are both  disjoint weak$^*$ null $L$\,-sequence in $(E^{\prime})^{+}$. From the inequality
	$$\|{T^{\prime}f_{n}}\|\leq{\|{T^{\prime}(f_{n}^{+})}\|+{\|{T^{\prime}(f_{n}^{-})}\|}}$$
	it follows that $\|{T^{\prime}f_{n}}\|\to{0}$.
\end{proof}
\begin{remark}\label{remark disjointly}
i) (1)$\Leftrightarrow$(3) establishes our conjecture that $T:X\to{E}$ is an almost $(L)$ limited operator if and only if $T$ is  disjointly weakly compact. Also, $T:X\to{E}$ is an almost $(L)$ limited operator if and only if $T^{\prime}$ is order weakly compact. Here, $E$ is not assumed to be $\sigma$-Dedekind complete.
\par ii)  (2)$\Leftrightarrow$(4) was recently proved by Elbour in \cite{El}.
\par iii) Let $T:X\to{E}$ be a disjointly weakly compact operator and let $S:E\to{F}$ be an arbitrary order bounded operator. Then the composition operator $ST:X\to{F}$ is also disjointly weakly compact. This follows easily from the equivalence (1)$\Leftrightarrow$(4) in Theorem \ref{Theorem 3.2}.
\end{remark}

Note that $E^{\prime}$ has order continuous norm if and only if the closed unit ball $B_E$ is disjointly weakly compact, that is,  the identity operator $id_{E}$ on $E$ is disjointly weakly compact.  As an immediate consequence of  Theorem \ref{Theorem 3.2}, the following result characterizes  dual Banach lattices with order continuous norms.
\begin{corollary}\label{Corollary 3.3}
	For a Banach lattice $E$ the following statements are equivalent.
	\begin{enumerate}
		\item The norm of $E^{\prime}$ is order continuous, i.e, $B_E$ is disjointly weakly compact.
		\item $\|{f_{n}}\|\to{0}$ for each order bounded sequence $(f_{n})\subset{(E^{\prime})^{+}}$ satisfying $f_{n} \stackrel{w^*}{\longrightarrow}0$.
		\item $\|{f_{n}}\|\to{0}$ for each disjoint $L$-sequence $(f_{n})\subset{E^{\prime}}$.
		\item $\|{f_{n}}\|\to{0}$ for each weak$^*$ null $L$\,-sequence $(f_{n})\subset{(E^{\prime})^{+}}$.
	\end{enumerate}	
\end{corollary}
\begin{remark}
Among the implications in Corollary \ref{Corollary 3.3},  (1)$\Leftrightarrow$(2) is also included in Theorem 3.1 of \cite{Wnuk}.
\end{remark}
\par The following result characterizes Banach lattices for which every  operator is disjointly weakly compact.
\begin{theorem}\label{Theorem 3.6}
	For a Banach lattice $E$ the following statements are equivalent.
	\begin{enumerate}
		\item Each continuous operator $T:X\to{E}$ from an arbitrary Banach space $X$ into $E$ is disjointly weakly compact.
		\item Each continuous operator $T:{\ell_{1}}\to{E}$ is disjointly weakly compact.
		\item The norm of $E^{\prime}$ is order continuous.
	\end{enumerate}	
\end{theorem}
\begin{proof}
	(1)$\Rightarrow$(2) and (3)$\Rightarrow$(1) Obvious.\par
	(2)$\Rightarrow$(3)
	we assume by way of contradiction that the norm of $E^{\prime}$ is not order continuous. Then from Corollary \ref{Corollary 3.3} it follows that there exists a disjoint $L$\,-sequence $(f_{n})$ of $E^{\prime}$ which is not norm null.
	By passing to a subsequence if necessary, we may assume that there exists some $\epsilon>{0}$ and $(x_{n})\subset{(B_{E})^{+}}$ such that $|{f_{n}(x_{n})}|\geq{\epsilon}$ for all $n$. Now consider the operators $T:{\ell_{1}}\to{E}$ defined by
	$$T(({\lambda}_{i})_{i=1}^{\infty})=\sum_{i=1}^{\infty}{\lambda_{i}}{x_{i}}$$ for each $(\lambda_{i})_{i=1}^{\infty}\in{\ell_{1}}$.
	Clearly, $T$ is well-defined, and $T^{\prime}f=(f(x_{i}))_{i=1}^{\infty}$ for all $f\in E^{\prime}$. We claim that $T$ is not disjointly weakly compact. Indeed, for the disjoint $L$\,-sequence $(f_{n})\subset E^{\prime}$, we have
	$$\|{T^{\prime}f_{n}}\|=\sup\limits_{i}{|{f_{n}(x_{i})}|}\geq{|{{f_{n}}(x_{n})}|}\geq{\epsilon}\quad \textmd{for all}\,\, n.$$This leads to a  contradiction.	
\end{proof}

\par A disjointly weakly compact operator between Banach lattices need not be order weakly compact, and vice versa. For instance, the identity operator $id_{\ell_{1}}:{\ell_{1}}\to{\ell_{1}}$ is order weakly compact, but not disjointly weakly compact. On the other hand, the identity operator $id_{\ell_{\infty}}:{\ell_{\infty}}\to{\ell_{\infty}}$ is disjointly weakly compact, but not order weakly compact.

\begin{theorem}\label{Theorem 4.2}
	Let $E$ and $F$ be two Banach lattices. Then the following assertions are equivalent.	
	\begin{enumerate}
		\item Every order bounded operator $T:E\to{F}$ is disjointly weakly compact.
		\item Every order bounded order weakly compact operator $T:E\to{F}$ is disjointly weakly compact.
		\item Every positive almost Dunford-Pettis operator $T:E\to{F}$ is disjointly weakly compact.
		\item One of the following conditions is valid:
		\begin{enumerate}
			\item the norm of $E^{\prime}$ is order continuous.	
			\item the norm of $F^{\prime}$ is order continuous.		
		\end{enumerate}	
	\end{enumerate}	
\end{theorem}
\begin{proof}
	(1)$\Rightarrow$(2)$\Rightarrow$(3) Obvious.\par
	(3)$\Rightarrow$(4) It suffices to establish that if the norm of $E^{\prime}$ is not order continuous, then $F^{\prime}$ has an order continuous norm.\par
	Assume that the norm of $E^{\prime}$ is not order continuous. Then  $E$ contains a sublattice isomorphic to $\ell_{1}$ (cf. \cite[Theorem 2.4.14 ]{Meyer}). It follows from \cite[Proposition 2.3.11]{Meyer} that there exists a positive projection $P:E\to{\ell_{1}}$. Since the operator $P$ is surjective, there exist $\delta>{0}$ such that $\delta\cdot{B_{\ell_{1}}}\subset{P(B_{E})}$. On the other hand, for every disjoint sequence $(y_{n})$ in $B_{F}^{+}$, the operator $Q:\ell_{1}\to{F}$ is defined by
	\begin{center}
		$ Q(\alpha_{n})=\sum_{n=1}^{\infty}{\frac{\alpha_{n}}{\delta}}{y_{n}}$ for each $(\alpha_{n})\in{\ell_{1}}$.
	\end{center}
	Clearly, $Q$ is well defined and positive.	
	Now, we consider the composed operator $T=Q\circ{P}: E\to{\ell_{1}}\to{F}$. Clearly,  $T$ is positive and almost Dunford-Pettis because $\ell_{1}$ has the  Schur property. Thus by hypothesis $T$ is a disjointly weakly compact operator. Also,  we have
	\begin{center}
		$y_{n}=\delta\cdot Q(e_{n})\in\delta\cdot Q(B_{\ell_{1}})$ and $\delta\cdot Q(B_{\ell_{1}})\subset Sol(T(B_{E}))$
	\end{center}
	(where $(e_{n})_{n=1}^{\infty}$ is the  standard vector basis of $\ell_{1}$). Therefore  ${y_n}\stackrel{w}{\rightarrow}0$ since $T(B_{E}$ is disjointly weakly compact.
	This implies that the norm of $F^{\prime}$ is order continuous (see, e.g., \cite[Theorem 2.4.14]{Meyer}).\par
	(4; a)$\Rightarrow$(1) Let $T:E\to{F}$ be an order bounded operator. Then $T^{\prime}:F^{\prime}\to{E^{\prime}}$ is also order bounded. If $E^{\prime}$ has order continuous norm, then $T^{\prime}$ is order weakly compact. By Theorem \ref{Theorem 3.2}, $T$ is disjointly weakly compact.\par
	(4; b)$\Rightarrow$(1) Follows directly from Theorem \ref{Theorem 3.6}.
\end{proof}
 By definition an operator between Banach spaces is called weakly precompact if it sends norm bounded sets to weakly precompact ones. It should be noted that every weakly precompact operator from a Banach space to a Banach lattice is disjointly weakly compact. For the converse we have the following.
\begin{theorem}\label{Theorem 4.5}
	Let $E$ and $F$ be Banach lattices. If every positive disjointly weakly compact operator from $E$ into $F$ is weakly precompact, then one of the following assertions is valid:
	\begin{enumerate}
       	\item The norm of $E^{\prime}$ is order continuous.

		\item Every order interval in $F$ is weakly precompact.
	\end{enumerate}	
\end{theorem}
\begin{proof}
	Assume by way of contradiction that neither (1) nor (2) holds. It follows from Theorem 3.6 of \cite{Xiang} that there exist two positive operators $S, T:E\to{F}$ such that $0\le{S}\le{T}$ holds and $T$ is weakly precompact while $S$ is not weakly precompact. The weak precompactness of $T$  implies that $T$ is disjointly weakly compact. Therefore $S$ is likewise disjointly weakly compact. By our hypothesis $S$ is weakly precompact. This is impossible.
\end{proof}
\begin{remark}
One may ask whether a positive disjointly weakly compact operator is weakly precompact if one of the two assertions in Theorem \ref{Theorem 4.5} is valid. The first assertion can not guarantee the weak precompactness of  disjointly weakly compact operators. For instance,  the identity operator $id_{\infty}:\ell_{\infty}\to{\ell_{\infty}}$ is not weakly precompact while it is disjointly weakly compact. About the second assertion in Theorem \ref{Theorem 4.5} we can say more. The following result improves Corollary 3.4.14 of \cite{Meyer}.
\end{remark}
\begin{corollary}\label{Corolloary 4.6}
	Let $F$ be a Banach lattice. Then the following statements are equivalent.
	\begin{enumerate}
		\item Every order interval in $F$ is weakly precompact.

		\item Every disjointly weakly compact set of $F$ is weakly precompact.
        \item Every disjointly weakly compact operator from an arbitrary Banach space $X$ to $F$ is weakly precompact.
		\item Every positive disjointly weakly compact operator from $\ell_{1}$ to $F$ is weakly precompact.
	\end{enumerate}	
\end{corollary}
\begin{proof} 	
(1)$\Leftrightarrow$(2) is included in Theorem 2.4 of \cite{Xiang}.\par
(2)$\Rightarrow$(3)$\Rightarrow$(4) Obvious.\par
(4)$\Rightarrow$(1) Follows from Theorem \ref{Theorem 4.5} and the fact that the norm of $\ell_{1}^{\prime}=\ell_{\infty}$ is not order continuous.
\end{proof}

It is well known that a Banach lattice has order continuous norm if and only if its order intervals are all weakly compact. In the following lemma we characterize the  order continuity of norms by means of weak precompactness.

\begin{lemma}\label{order continuity}
	The norm of a Banach lattice $E$ is order continuous if and only if $E$ is $\sigma$-Dedekind complete and each order interval of $E$ is weakly precompact. In particular, $E^{\prime}$ has order continuous norm if and only if  the order intervals in $E^{\prime}$ are all weakly precompact.	
\end{lemma}
\begin{proof}
	It suffices to prove that  $E$ has  order continuous norm if $E$ is $\sigma$-Dedekind complete and each order interval of $E$ is weakly precompact.	To this end, we assume by way of contradiction that $E$ does not have order continuous norm. Then $\ell_{\infty}$ is lattice embeddable in $E$. Without loss of generality we can assume that $\ell_{\infty}$ is a closed sublattice of $E$. Since every order bounded subset of $E$ is weakly precompact, it follows that $B_{\ell_{\infty}}$ is also weakly precompact. This is impossible.
\end{proof}

\begin{theorem}
	For $\sigma$-Dedekind complete Banach lattices $E$ and $F$ the following conditions are equivalent.
 \begin{enumerate}
\item  Each  disjointly weakly compact operator $T:E\to{F}$  is weakly precompact.
\item One of the following assertions is valid:
\begin{enumerate}
		\item The closed unit ball $B_{E}$ is weakly precompact.
		\item The norm of $F$ is order continuous.
\end{enumerate}	
\end{enumerate}
\end{theorem}
\begin{proof} (1)$\Rightarrow$(2) It suffices to show that if  the norm of $F$ is not order continuous, then $B_{E}$ is weakly precompact. Firstly, from Lemma \ref{order continuity} and Theorem \ref{Theorem 4.5} it follows that the norm of $E^{\prime}$ is order continuous. Secondly, by our hypothesis and Proposition 3.4.15 of \cite{Meyer} every disjointly weakly compact operator $T:E\to{F}$ is also order weakly compact. Since disjointly weakly compact operators are just  those operators whose adjoint operators are order weakly compact,  it follows from  \cite[Theorem 2.8]{JHm} that $E$ has order continuous norm.
Finally, we can see that the unit ball $B_{E}$ is weakly precompact because the norms on $E$ and $E^{\prime}$ are order continuous ( cf. \cite[Corollary 2.5.10]{Meyer}.\par
(2; a)$\Rightarrow$(1) Obvious.\par
(2; b)$\Rightarrow$(1) Let $T:E\to{F}$ be  a disjointly weakly compact operator. Since $F$ has order continuous norm, then every order interval in $F$ is weakly precompact. It follows from \cite[Theorem 2.4]{Xiang} that the disjointly weakly compact set $T(B_{X})$ is weakly precompact, that is, $T$ is weakly precompact.
\end{proof}

\par Recall that a Banach lattice $E$ has the $\sigma$-Levi property whenever $$0\leq x_{n}\uparrow\,\,  \textrm{and}\,\, \sup_{n}\{\|x_{n}\|\}<\infty \,\Rightarrow \sup_{n}\{x_{n}\} \, \textrm{exists in \textit{E}}.$$ Every $KB$-space and every $\sigma$-Dedekind complete $AM$-space with unit has the $\sigma$-Levi property. A Banach lattice is a $KB$-space if and only if it has both order continuous norm and the $\sigma$-Levi property.
Next we characterize Banach lattices for which each positive disjointly weakly compact operator is weakly compact.
\begin{theorem}\label{Theorem 4.7}
	Let $E$ and $F$ be two Banach lattices such that $F$ has the $\sigma$-Levi property . Then the following assertions are equivalent:
	\begin{enumerate}
		\item Every disjointly weakly compact operator $T:E\to{F}$ is weakly compact.
		\item Every positive disjointly weakly compact operator $T:E\to{F}$ is weakly compact.
		\item One of the following conditions is valid:
		\begin{enumerate}
			\item  $E$ is reflexive.
			\item  $F$ is a $KB$-space.
		\end{enumerate}
	\end{enumerate}
\end{theorem}
\begin{proof} (1)$\Rightarrow$(2) Obvious.\par
	(2)$\Rightarrow$(3) It suffices to establish that if $E$ is not reflexive, then $F$ is a $KB$-space. Since $F$ has the $\sigma$-Levi property, it remains to prove that the norm of $F$ is order continuous. \par
	 Since every semi-compact operator is disjointly weakly compact (see Theorem \ref{Theorem 3.2} and \cite[Proposition 3.6.18]{Meyer} for details), then, by our hypothesis, every positive semi-compact operator is weakly compact. Since every Banach lattice with the $\sigma$-Levi property is $\sigma$-Dedekind complete, it follows from Theorem 2.2 of \cite{Aqzzouz5} that $F$ has order continuous norm.\par	
	(3; a)$\Rightarrow$(1) Obvious.\par
	(3; b)$\Rightarrow$(1) Follows from Proposition \ref{Proposition 2.5}.
\end{proof}	
\begin{remark}\label{Remark 4.8} The assumption  ``$F$ has the $\sigma$-Levi property" is essential in Theorem \ref{Theorem 4.7}. For instance, $F=c_{0}$ has order continuous norm, but  does not have the  $\sigma$-Levi property. Note that  each bounded linear operator $T:\ell_{\infty}\to{c_{0}}$ is weakly compact (because $\ell_{\infty}$ is a Grothendieck space). However, $\ell_{\infty}$ is not reflexive and $c_{0}$ is not a $KB$-space.
\end{remark}
\par Let us recall that a Banach lattice $E$ has
\begin{itemize}
	\item[$\raisebox{0mm}{------}$] the \textit{dual positive Schur property} whenever $0\leq f_{n}\stackrel{w^*}{\longrightarrow}0$ implies  $\|f_{n}\|\to{0}$ in $E^{\prime}$ (see \cite{AEW, Wnuk2}).

   \item[$\raisebox{0mm}{------}$] the \textit{positive Grothendieck property} if every weak$^*$ null sequence in $E^{\prime}$ with positive terms
is weakly null.
\end{itemize}
It is very clear that a Banach lattice $E$ has the dual positive Schur property if and only if $E$ has the positive Grothendieck property and $E^{\prime}$ has the  positive Schur property (see \cite{Wnuk2}). Next we come to $L$\,-weak compactness of  disjointly weakly compact operators. As a consequence of  \cite[Theorem 3.4]{ZLC} and Corollary \ref{psp},  the following result is immediate.
\begin{proposition}\label{Proposition 4.10}
	Let $E$ be a Banach lattice. Then the following statements are equivalent.
	\begin{enumerate}
       \item Every disjointly weakly compact operator $T:X\to E$ is $L$\,-weakly compact.
		\item Every positive disjointly weakly compact operator $T:\ell_{1}\to E$ is $L$\,-weakly compact.
        \item Every positive weakly compact operator $T:\ell_{1}\to E$ is $L$\,-weakly compact.
		\item $E$ has the positive Schur property.
	\end{enumerate}		
\end{proposition}
Recall from \cite[Theorem 3.1]{BAJ} that a Banach lattice $E$ is finite dimensional if and only if the closed unit ball $B_{E}$ of $E$ is $L$\,-weakly compact. On the other hand, a Banach lattice $E$ is finite dimensional if and only if it is reflexive and it has the positive Schur property (see the proof of \cite[Theorem 2]{CNa} for details). It follows from Proposition 2.1\,(b) of \cite{Wnuk2} that a Banach lattice $E$ is finite dimensional if and only if it has order continuous norm and  the dual positive Schur property.
\begin{theorem}\label{Theorem 4.11}
	For a Banach lattice $E$  the following statements are equivalent.
	\begin{enumerate}
		\item $E$ has the positive Schur property.
        \item Every positive disjointly weakly compact operator on $E$ is $L$\,-weakly compact.
	\end{enumerate}		
\end{theorem}	
\begin{proof}
(1)$\Rightarrow$(2) Obvious.\par
(2)$\Rightarrow$(1) We assume by way of contradiction that $E$ does not have the positive Schur property.
Firstly, by our hypothesis it is clear that every positive weakly compact operator on $E$ is $L$\,-weakly compact. It follows from \cite[Theorem 3.5]{Aqzzouz2} that $E^{\prime}$ has the  positive Schur property and $E$ has an order continuous norm.\par
 Secondly, we claim that $E$ is not a $KB$-space. Otherwise, $E$ is reflexive, and hence the identity operator $id_{E}$ is  (disjointly) weakly compact. By our hypothesis again $id_{E}$ is  $L$\,-weakly compact, that is, $B_E$ is $L$\,-weakly compact. Then it follows from \cite[Theorem 3.1]{BAJ} that $E$ is finite dimensional. This is impossible. Let $i:c_{0}\to{E}$ be the lattice embedding.\par
 Thirdly, we claim that $E$ has the positive Grothendieck property. Otherwise, it follows from Proposition 2.12 of \cite{Wnuk2} that there exists a positive surjection $T:E\to{c_{0}}$. By the Open Mapping Theorem  we can find some $\delta>{0}$ such that $\delta\cdot{B_{c_{0}}}\subset{T(B_{E})}$. Clearly, $T$ is disjointly weakly compact since $(c_{0})^{\prime}={\ell_{1}}$ has order continuous norm.  Then, $i\circ T: E\stackrel{T}{\rightarrow}c_{0}\stackrel{i}{\rightarrow}E$ is also disjointly weakly compact (see Remark \ref{remark disjointly} iii)). By our hypothesis, $i\circ T$ is $L$\,-weakly compact. This implies that $\delta\cdot i({B_{c_{0}}})\, \subset{i\circ T(B_{E})}$ is $L$\,-weakly compact. This is impossible. Therefore $E$ has the positive Grothendieck property, and hence $E$ has the dual positive Schur property since $E^{\prime}$ has the positive Schur property.\par
 To sum up, we have shown that if we assume that $E$ does not have the positive Schur property, then $E$ has an order continuous norm and  the  dual positive Schur property. It follows from Proposition 2.1 (b) of \cite{Wnuk2} that $E$ is finite dimensional. This is impossible. Therefore, $E$ has the positive Schur property.
\end{proof}
\par We conclude with two results concerning $M$-weak compactness of disjointly weakly compact operators. Let us recall that
\begin{itemize}
	\item[$\raisebox{0mm}{------}$] a continuous operator $T:E\to{X}$ from a Banach lattice into a Banach space is said to be \textit{$M$-weakly compact} whenever $\lim\|Tx_n\|=0$ holds for every norm  bounded disjoint sequence $(x_n)\subset E$.

   \item[$\raisebox{0mm}{------}$] a continuous operator $T:X\to{E}$ from a Banach space into a Banach lattice is \textit{semi-compact} whenever for each $\epsilon>0$ there exists some $u\in E^{+}$ such that $T(B_X)\subset [-u,u]+\epsilon B_B$.
\end{itemize}
 Every $M$-weakly compact operator is weakly compact. Positive $M$-weakly compact operators between Banach lattices are semi-compact.  An operator $T:E\to{X}$ from an $M$-space into a Banach space is $M$-weakly compact if and only if it is weakly compact (see, e.g., pp. 321-322, 331-333 in \cite{AB} for details). It should be noted that every semi-compact operator is disjointly weakly compact.
\begin{theorem}\label{Theorem 4.12}
	For a Banach lattice $E$ the following statements are equivalent.
	\begin{enumerate}
		\item $E$ is finite dimensional.
		\item Every positive disjointly weakly compact operator $T:E\to\ell_{\infty}$ is $M$-weakly compact.
	\end{enumerate}
\end{theorem}
\begin{proof}
(1)$\Rightarrow$(2) Obvious.\par
(2)$\Rightarrow$(1) We assume by way of contradiction that $E$ is infinite dimensional. Then there exists a positive disjoint sequence $(x_{n})$ of $E^{+}$ such that $\|x_{n}\|=1$ for all $n$. Therefore, there exists   a  disjoint sequence $(g_{n})$ of $(E^{\prime})^+$  such that
\begin{center}
 $\|g_{n}\|=1$ and $g_{n}(x_{n})=1$  for all $n$,    \quad $g_{n}(x_{m})=0$  for $n\neq{m}$.
\end{center}
We define the positive operator $T:E\to\ell_{\infty}$ by
$T(x)=(g_{1}(x),g_{2}(x),\cdots)$
for all $x\in{E}$. Clearly, $T$  is a  disjointly weakly compact positive operator. However, $T$ is not $M$-weakly compact. Indeed, for the disjoint sequence $(x_{n})$ of $B_{E}^{+}$, we have
$$\begin{aligned}
	\|T(x_{n})\|=\sup_{i}|g_{i}(x_{n})|={|g_{n}(x_{n})|}=1.	
	\end{aligned}$$
This contradicts with our hypothesis. Hence, $E$ is a finite dimensional Banach lattice.
\end{proof}
\begin{theorem}\label{M-weak}
Let $E$ be a $\sigma$-Dedekind complete Banach lattice. The following statements are equivalent.
\begin{enumerate}
		\item $E$ is finite dimensional.
        \item Every positive disjointly weakly compact operator on  $E$ is $M$-weakly compact.
 \end{enumerate}
\end{theorem}
\begin{proof}
(1)$\Rightarrow$(2) Obvious.\par
(2)$\Rightarrow$(1) If every positive disjointly weakly compact operator on  $E$ is $M$-weakly compact, then every positive semi-compact operator on  $E$ is $M$-weakly compact since semi-compact operators are disjointly weakly compact. It follows from \cite[Theorem 2]{CNa} that at least one of the following two assertions holds:\par
 i) $E$ is finite dimensional.\par
 ii) The norms of  $E^{\prime}$ and $E$ are order continuous. \par
 If $E$ is finite dimensional, it is  done. If the norms of  $E^{\prime}$ and $E$ are order continuous, then, by a result of P. G. Dodds and D. H. Fremlin,  the class of positive $M$-weakly compact operators on $E$ and the class of positive $L$\,-weakly compact operators coincide (see, e.g., \cite[Theorem 5.67]{AB}).  This implies that  every positive disjointly weakly compact operator on  $E$ is $L$\,-weakly compact. From Theorem \ref{Theorem 4.11} it follows that $E$ has the positive Schur property. Again, the order continuity of the norm on $E^{\prime}$ implies that $B_E$ is disjointly weakly compact. Then, by Corollary \ref{psp} $B_E$ is $L$\,-weakly compact, i.e., $E$ is finite dimensional.
\end{proof}


\vskip 5.5mm
\begin{thebibliography}{99}
	


	
\bibitem{AB}C.D. Aliprantis and O.~Burkinshaw, \textit{Positive Operators} (Reprint of the 1985 original), Springer, Dordrecht, 2006.
\bibitem{Andrews} K. Andrews, Dunford-Pettis sets in the space of Bochner integrable functions, \textit{Math. Ann.} \textbf{241} (1979),  35-41.
\bibitem{Aqzzouz1}B.~Aqzzouz and K.~Bouras, $(L)$ sets and almost $(L)$ sets in Banach lattices,\textit{ Quaest. Math.} \textbf{36} (2013), 107-118.
\bibitem{JHm}
B.~Aqzzouz and J.~Hmichane, The duality problem for the class of order weakly compact operators, \textit{Glasg. Math. J.} \textbf{51} (2009), 101-108.
\bibitem{BAJ}
B.~Aqzzouz, A.~Elbour and J.~Hmichane, Some properties of the class of positive Dunford-Pettis operators,
\textit{J. Math. Anal. Appl.} \textbf{354} (2009), 295-300.
\bibitem{Aqzzouz2}
B.~Aqzzouz and A.~Elbour, The Schur property of Banach lattices and the compactness of weakly compact operators, \textit{Math. Proc. R. Ir. Acad.} \textbf{110A} (2010), 1-11.

\bibitem{AEW}B. Aqzzouz, A. Elbour and A.W. Wickstead, Positive almost Dunford-Pettis operators and their duality,
\textit{Positivity} \textbf{15} (2011), 185-197.
\bibitem{Aqzzouz5} B.~Aqzzouz and A.~Elbour, On the weak compactness of semi-compact operators, \textit{Math. Rep. (Bucur.)} \textbf{16} (2014), 1-6.

 \bibitem{BFV} J. Borwein, M. Fabian and J. Vanderwerff, Characterizations of Banach spaces via convex and other locally Lipschitz functions, \textit{Acta Math. Vietnam} \textbf{22} (1997), 53-69.

\bibitem{Bour}
K.~Bouras, Almost Dunford-Pettis sets in Banach lattices, \textit{Rend. Circ. Mat. Palermo } \textbf{62} (2013), 227-236.
\bibitem{BD}J. Bourgain and J. Diestel, Limited operators and strict cosingularity, \textit{Math. Nachr.} \textbf{119} (1984), 55-58.
\bibitem{CGL}H. Carri\'{o}n, P. Galindo and M.L. Louren\c{c}o, A stronger Dunford-Pettis property, \textit{Studia Math.} \textbf{184} (2008), 205--216.
\bibitem{CCJ}J.X. Chen, Z.L. Chen and G.X. Ji, Almost limited sets in Banach lattices, \textit{J. Math. Anal. Appl.} \textbf{412} (2014), 547-553.

\bibitem{ZLC}
Z.L. Chen and  A.W. Wickstead,  $L$-weakly and $M$-weakly compact operators, \textit{Indag. Math. (N.S.)} \textbf{10} (1999), 321-336.
\bibitem{CNa}
N.~Cheng, Z.L. Chen and Y.~Feng, $L$- and $M$-weak compactness of positive semi-compact operators, \textit{Rend. Circ. Mat. Palermo (2)}  \textbf{59} (2010), 101-105.
\bibitem{El} A. Elbour, Some characterizations of almost limited operators, \textit{Positivity} \textbf{ 21} (2017), 865-874.
\bibitem{FMM} K. El Fahri, N. Machrafi and M. Moussa, Banach lattices with the positive Dunford-Pettis relatively compact property, \textit{Extracta Math.} \textbf{30} (2015), 161-179.

\bibitem{Em}G. Emmanuele, A dual characterization of Banach spaces not containing $\ell^{1}$, \textit{Bull. Pol. Acad. Sci. Math.} \textbf{34} (1986), 155-160.

\bibitem{Ghenciu} I.~Ghenciu, A note on weak reciprocal Dunford-Pettis sets, \textit{Acta Math. Hungar.} \textbf{152} (2017), 453-463.
\bibitem{Leung} D.H. Leung, On the weak Dunford-Pettis property, \textit{Arch. Math.}\textbf{ 52} (1989), 363-364.
\bibitem{MEM} N. Machrafi, A. Elbour and M. Moussa, Some characterizations of almost limited sets and
applications. arXiv:1312.2770v1

\bibitem{Meyer}P.~Meyer-Nieberg, \textit{Banach Lattices}, Universitext, Springer-Verlag, Berlin, 1991.
\bibitem{Oughajji}
F.Z. Oughajji, K.~El Fahri and M.~Moussa, On the class of almost $(L)$ limited operators, \textit{Acta Sci. Math. (Szeged)} \textbf{87} (2021), 207-218.
\bibitem{Sa} J.A. Sanchez, \textit{Operators on Banach Lattices} (Spanish), Ph.D. Thesis, Complutense University, Madrid, 1985.

 \bibitem{WC} Y. Wen and J.X. Chen, Characterizations of Banach spaces with relatively compact Dunford-Pettis sets, \textit{Adv. Math. (China) } \textbf{45} (2016), 122-132.

\bibitem{Wnuk}
W.~Wnuk, \textit{Banach lattices with order continuous norms,} Polish Scientific Publishers PWN, Warsaw (1999).
\bibitem{Wnuk2}
W.~Wnuk, On the dual positive Schur property in Banach lattices, \textit{Positivity} \textbf{17} (2013), 759-773.
\bibitem{Xiang}
B.~Xiang, J.X. Chen and L.~Li, Weak precompactness in Banach lattices, \textit{Positivity} \textbf{26} (1) (2022), 4. https://doi.org/10.1007/s11117-022-00873-6



	% etc
\end{thebibliography}


\end{document}